\documentclass[11pt]{amsart}

\usepackage{epsfig}
\usepackage{graphicx}
\usepackage[caption=false,font=footnotesize]{subfig}
\usepackage{verbatim}
\usepackage{amssymb}
\usepackage{algorithm,algorithmic}
\usepackage{tikz}

\newcommand{\bS}{\mathbf{S}}
\newcommand{\rmn}{\mathrm{n}}
\newcommand{\polylog}{\textup{polylog}}
\newcommand{\OO}{\Omega}
\newcommand{\oo}{\omega}

\renewcommand{\implies}{\quad \Longrightarrow \quad}

\newcommand{\bigo}{\textup{O}}

\renewcommand{\phi}{\varphi}
\newcommand{\eps}{\varepsilon}

\newcommand{\E}{\mathbb{E}}

\newcommand{\Var}{\mathrm{Var}}

\newcommand{\e}{\textup{e}}
\newcommand{\wh}[1]{\widehat{#1}}

\newcommand{\poly}{\textup{poly}}
\DeclareMathOperator*{\Arg}{Arg}
\newcommand{\bm}[1]{\mathbf{#1}}

\newcommand{\wt}{\widetilde}

\newcommand{\eproof}{\hfill\rule{2.2mm}{3.0mm}}

\newcommand{\Proof}{\noindent {\bf Proof.~~}}

\newcommand{\R}{{\mathbb R}}
\newcommand{\Z}{{\mathbb Z}}

\newcommand{\C}{{\mathbb C}}

\newcommand{\wmod}[1]{\mbox{~(mod~$#1$)}}
\newcommand{\wmodc}{\wmod{[-\frac{1}{2}, \frac{1}{2})}}
\renewcommand{\eqref}[1]{(\ref{#1})}

\newcommand{\shsp}{\hspace{1em}}
\newcommand{\mhsp}{\hspace{2em}}

\newtheorem{prop}{Proposition}[section]
\newtheorem{lem}[prop]{Lemma}

\newtheorem{coro}[prop]{Corollary}
\newtheorem{theo}[prop]{Theorem}

\begin{document}
\baselineskip 18pt

\title{A Multiscale Sub-linear Time Fourier Algorithm for Noisy Data}

\author{Andrew~Christlieb}
\address{Department of Mathematics  \\ Michigan State University \\
East Lansing, MI 48824, USA.}
\email{christlieb@math.msu.edu}

\author{David~Lawlor}
\address{Statistical and Applied Mathematical Sciences Institute \\ Department of Mathematics  \\ Duke University \\
Durham, NC 27708, USA.}
\email{djl@math.duke.edu}

\author{Yang~Wang}
\address{Department of Mathematics  \\ Michigan State University \\
East Lansing, MI 48824, USA.}
\email{ywang@math.msu.edu}
%%%%%%%%%%%%%%%%%%%%%%%%

\maketitle

\begin{abstract}
We extend the recent sparse Fourier transform algorithm of \cite{lawlor2013adaptive} to the noisy setting, in which a signal of bandwidth $N$ is given as a superposition of $k\ll N$ frequencies and additive noise. We present two such extensions, the second of which exhibits a novel form of error-correction in its frequency estimation not unlike that of the $\beta$-encoders in analog-to-digital conversion \cite{daubechies2006d}. The algorithm runs in time $\bigo(k \log(k) \log(N/k))$ on average, provided the noise is not overwhelming. The error-correction property allows the algorithm to outperform FFTW~\cite{frigo2005design}, a highly optimized software package for computing the full discrete Fourier transform, over a wide range of sparsity and noise values, and is to the best of our knowledge novel in the sparse Fourier transform context. 
\end{abstract}

%%%%%%%%%%%%%%%%%%%%%%
\section{Introduction} \label{sec:intro}
\setcounter{equation}{0}

The Fast Fourier Transform (FFT) \cite{cooley1965algorithm} is a fundamental numerical algorithm whose importance in a wide variety of applications cannot be overstated. The FFT reduces the runtime complexity of calculating the discrete Fourier transform (DFT) of a length $N$ array from the naive $\bigo(N^2)$\footnote{We write $f=\bigo(g)$ to indicate that $f(x) \le cg(x)$ for some positive constant $c$ and all sufficiently large $x$.} to $\bigo(N\log(N))$. At the time of its introduction in the mid-1960s, it dramatically increased the size of problems that a typical computer could handle. Over the past fifty years the typical size of data sets has grown by orders of magnitude, and in certain application areas (e.g.~ultra-wideband radar) the computation of the full FFT is no longer tractable on commodity hardware. In this and other instances, however, it is known \emph{a priori} that the signals of interest have small frequency support; that is, their Fourier transforms are \emph{sparse}. This problem has received attention from a number of research communities over the past decade, who have shown that it is possible to significantly outperform the FFT in both runtime and sampling requirements when the number of significant Fourier modes $k$ is much less than the nominal bandwidth $N$.

The earliest work to specifically address the sparse Fourier transform problem was \cite{gilbert2002near}, which gave a randomized algorithm with runtime and sampling complexity $\bigo(k^2\,\polylog(N))$.\footnote{We write $f=\polylog(g)$ to indicate that $f=\bigo(\log^c(g))$ for some unspecified constant $c$.} This was later improved to $\bigo(k\,\polylog(N))$ \cite{gilbert2005improved} through the use of unequally-spaced FFTs \cite{anderson1996rapid}. For a given failure probability $\delta$ and accuracy parameter $\eps$, the algorithm returns a $k$-term approximation $\hat{y}$ to the DFT of the input $\hat{x}$ such that with probability $1-\delta$ it holds that
\begin{equation}
\label{eq:error-ggi+}
\|\hat{x}-\hat{y}\|_2^2 \le (1+\eps)\|\hat{x}-\hat{x}_k\|_2^2.
\end{equation}
Here $\hat{x}_k$ is the best $k$-term approximation to $\hat{x}$ and $\|\cdot\|_2$ is the discrete $\ell_2$ norm. A separate group
of authors \cite{hassanieh2012simple} has developed a modified version of this algorithm with runtime
$\bigo(\log(N)\sqrt{Nk\log(N)})$. While the dependence on $N$ is sub-optimal asymptotically,
in practice this algorithm is significantly faster than either \cite{gilbert2002near}
or \cite{gilbert2005improved}. The same authors presented an improved algorithm with runtime $\bigo(k\log(N)\log(N/k))$ in \cite{hassanieh2012nearly} whose frequency identification prodecure is very similar to \cite{lawlor2013adaptive}, upon which the present work is based. However, the performance of \cite{hassanieh2012nearly} in the presence of noise has yet to be evaluated empirically. 

The algorithms described in the previous paragraph are all randomized, and so will fail on some small subset of potential inputs. Recognizing this as a potential detriment in failure-intolerant applications, two authors have independently given deterministic algorithms for the sparse Fourier transform problem. In \cite{akavia2010deterministic} an algorithm with runtime $\poly(k,\log(N))$\footnote{Here, $\poly(\cdot)$ indicates an unspecified polynomial in its arguments.} was given where the exponent on $k$ is at least six. This high dependence on  $k$ renders the algorithm infeasible in practice, and it has not been implemented. In \cite{iwen2010combinatorial}, the combinatorial properties of aliasing among frequencies were exploited to give an algorithm with runtime and sampling complexity $\bigo(k^2\,\polylog(N))$. While this represented a major improvement over the theoretical runtime complexity of \cite{akavia2010deterministic}, in practice it only outperformed the FFT for relatively modest values of the sparsity $k$.

Most recently the authors of \cite{lawlor2013adaptive} gave a deterministic algorithm with average-case sampling and runtime complexity $\bigo(k\log(N))$. The worst-case runtime
bounds are asymptotically of the same order as \cite{iwen2010combinatorial}, but over a representative
class of random signals it was shown to significantly outperform its deterministic and randomized competitors.
This was achieved by sampling the input at two sets of equispaced points slightly offset in time.
This time shift appears in the Fourier domain as a frequency modulation, which allows the authors
to both detect when aliasing has occurred and, for frequencies that are isolated (i.e.~not aliased),
to calculate the frequency value directly. While \cite{iwen2010combinatorial} also uses properties
of aliasing to reconstruct frequency values, it is not able to distinguish between aliased and
non-aliased terms until sufficiently many DFTs of coprime lengths have been computed, and so is
unable to perform any better in the average case than in the worst case. In the empirical evaluation
of \cite{lawlor2013adaptive} an improvement of over two orders of magnitude was observed over
\cite{gilbert2005improved} and \cite{iwen2010combinatorial}.

In this paper we extend the algorithm of \cite{lawlor2013adaptive} to noisy environments in
two distinct ways. The first of these, which is a minor modification of the noiseless algorithm, is based on a certain rounding of the frequency estimates and was previously reported in \cite{lawlor2013adaptive}. In this work we provide an improved algorithm and more detailed analysis of that earlier work. The second extension is the main result of this paper, a novel multiscale error-correcting
algorithm that utilizes offset time samples at geometrically spaced time shifts. This extension is
in essence a progressive frequency identification algorithm not unlike the $\beta$-encoders
for analog-to-digital conversion \cite{daubechies2006d}. The new algorithm gives excellent
performance in the noisy setting without significantly increasing the computational costs from
the noiseless case. For both extensions we provide detailed mathematical analysis as well as empirical evaluations. While both extensions work well in the noisy environment, the multiscale
algorithm achieves comparable accuracy at a significantly lower computational
cost.

The remainder of this paper is organized as follows. In Section~\ref{sec:prelim} we review the notation introduced in \cite{lawlor2013adaptive} that will be necessary in the sequel. We also describe our noise model, discuss some of the problems noisy signals present for the algorithm of \cite{lawlor2013adaptive}, and argue that in certain applications the $\ell_2$ error metric is inappropriate and should be replaced with a form of Earth Mover's Distance. We also describe the random signal model used in the empirical evaluations in Section~\ref{sec:empirical}. In Section~\ref{sec:minor} we give our first modified algorithm and analyze the dependence of the sampling rate on the noise level. In Section~\ref{sec:multi} we describe our multiscale frequency identification procedure, and in Section~\ref{sec:empirical} we provide an empirical evaluation of the accuracy and speed of both algorithms.  Finally in Section~\ref{sec:conc} we provide a brief conclusion.

%%%%%%%%%%%%%%%%%%%%%%%
\section{Preliminaries}
\label{sec:prelim}
\setcounter{equation}{0}

\subsection{Notation and brief review}
\label{subsec:notation}

In this section we introduce the notation that will be used in the remainder of this paper and briefly review the results in \cite{lawlor2013adaptive}. We denote by $\Z$ the set of integers, $\C$ the set of complex numbers, and we let $N$ be a fixed (large) natural number. We write $\lfloor x \rfloor$ to denote the largest integer less than or equal to $x$. All logarithms are in base two unless explicitly specified.

We consider frequency-sparse band-limited signals $S:[0,1)\to\C$ of the form
\begin{equation} \label{eq:signal}
     S(t) = \sum_{\oo\in\OO} a_\oo \e^{2\pi\textup{i}\omega t},
\end{equation}
where $\OO$ is a finite set of integers bounded in $[-N/2,N/2)$ and $0\neq a_\oo \in \C$ for each $\oo\in\OO$. For simplicity we shall extend $S(t)$ periodically to a function on the whole real line.
The  Fourier samples of $S$ are given by
\begin{equation} \label{eq:fourierseries}
   \widehat{S}(h) = \int_0^1 S(t) \e^{-2\pi\textup{i} h t} \textup{d}t, \; h \in \Z,
\end{equation}
so that for signals of the form~\eqref{eq:signal} we have $\widehat{S}(\oo) = a_\oo$ for $\oo\in\OO$
and $\wh S(h) = 0$ for all other $h\in \Z$. 

In practice we work with data of finite length. Given any finite  sequence ${\mathbf s}=(s_0, s_1, \dots, s_{p-1})$ of length $p$ its DFT is given by
\begin{equation} \label{eq:dft}
    \widehat{\bm{s}}[h] ~=~
        \sum_{j=0}^{p-1} s_j \e^{-2\pi\textup{i} jh/p}
        ~=~ \sum_{j=0}^{p-1} \bm{s}[j] W_p^{jh},
\end{equation}
where $h = 0, 1, \ldots, p-1$, $\bm{s}[j]:=s_j$ and $W_p:=\e^{-2\pi\textup{i}/p}$ is the primitive
$p$-th root of unity. The FFT \cite{cooley1965algorithm} allows the computation of $\wh{\bm{s}}$
in $\bigo(p\log p)$ steps.

All fast reconstruction algorithms apply the DFT to selected finite sample sets of $S(t)$,
and our work is no exception.  Let $p$ be a positive integer and $\eps>0$. The two sample
sets we use extensively are $\bS_p$ and $\bS_{p, \eps}$,
which are length $p$ samples of $S(t)$ given by
\begin{equation} \label{eq:sp}
    \bm{S}_p[j] = S\Bigl(\frac{j}{p}\Bigr), ~\shsp
    \bm{S}_{p,\eps}[j] = S\Bigl(\frac{j}{p}+ \eps\Bigr),~\shsp j = 0, 1, \ldots, p-1.
\end{equation}
For each $h$ let $\Lambda_{p,h} =\{\oo\in\OO: \oo\equiv h \wmod p\}$, where $\omega \equiv h \wmod{p}$ indicates that $\omega-h$ is divisible by $p$. It is a simple derivation to obtain
\begin{equation} \label{eq:sh}
      \wh\bS_{p}[h] = p \sum_{\oo\in\Lambda_{p,h}}  a_\oo,
        ~\shsp%   \mbox{~~and~~}
      \wh\bS_{p,\eps}[h] = p \sum_{\oo\in\Lambda_{p,h}}  a_\oo \e^{2\pi \textup{i}\eps \oo}.
\end{equation}
Let $\omega \wmod{p}$ indicate the remainder after division of $\omega$ by $p$. In the ideal scenario where all $\{\omega\wmod{p}:~ \oo\in\OO\}$ are distinct we have
\begin{equation}\label{eq:sph}
      \wh\bS_{p}[h] = \left\{\begin{array}{cl}
           pa_\oo& ~~h = \omega\wmod p \mbox{~for some~}\oo\in\OO, \\
                                0 & \mbox{~~otherwise},
      \end{array}\right.
\end{equation}
and similarly
\begin{equation}\label{eq:speh}
      \wh\bS_{p,\eps}[h] = \left\{\begin{array}{cl}
                             pa_\oo \e^{2\pi \textup{i}\eps \omega} & ~~h=\omega \wmod p
                             \mbox{~for some~$\oo\in\OO$}, \\
                                0 & \mbox{~~otherwise}.
      \end{array}\right.
\end{equation}
Thus, the nonzero elements of $\wh\bS_p[h]$ occur precisely at the locations $h = \omega \wmod{p}$ for some $\oo \in \OO$, and moreover for such $h$ we have $|\wh\bS_p[h]|=|\wh\bS_{p,\eps}[h]|$.
Furthermore for each $\oo\in\OO$ and $h = \omega \wmod p$ we have
$\frac{\wh\bS_{p,\eps}[h]}{\wh\bS_{p}[h]} = \e^{2\pi \textup{i} \eps\omega}$.
Hence
\begin{equation} \label{eq:omega1}
   2\pi\eps\omega \equiv \Arg\left(\frac{\widehat{\bS}_{p,\eps}[h]}{\widehat{\bS}_{p}[h]}\right)
   \wmod{2\pi},
\end{equation}
where $\Arg(z)$ denotes the phase angle of the complex number $z$ in $[-\pi, \pi)$.
Now assume that we have $|\eps| \leq \frac{1}{N}$. Then $\omega$ is completely determined by~\eqref{eq:omega1}, as there will be no wrap-around aliasing. Hence
\begin{equation} \label{eq:omega}
   \omega = \frac{1}{2\pi\eps} \Arg\left(\frac{\widehat{\bS}_{p,\eps}[h]}{\widehat{\bS}_{p}[h]}\right).
\end{equation}
The weight $a_{\oo}$ can be recovered via $a_{\oo} = \wh\bS_p[h]/p$.

\noindent
{\bf Remark.}~	
In fact, more generally, if we have an estimate of $\omega \in\OO$, say $|\omega| < \frac{L}{2}$, then by
taking $\eps \leq \frac{1}{L}$ the same reconstruction formula~\eqref{eq:omega} holds. We will use this
observation in Section~\ref{sec:multi} when we develop a multiscale frequency identification procedure for noisy signals.

Of course it is possible that not all $\{\omega \wmod{p}:~ \oo\in\OO\}$ are distinct. For an
$\oo\in\OO$ we say $\oo$ {\em has a collision modulo $p$}, or simply {\em has a collision}
when there is no ambiguity in the modulus $p$, if there is at least one other $\oo'\in\OO$ such that $\oo \equiv
\oo' \wmod p$. In \cite{lawlor2013adaptive} a criterion is developed to detect collisions in
the noiseless case. For $\oo\in\OO$ and $h=\omega \wmod p$, it is clear that a necessary
condition for no collision to occur is
\begin{equation}  \label{eq:non-collision}
\Bigl|\frac{\wh\bS_{p,\eps}[h]}{\wh\bS_{p}[h]}\Bigr| = |\e^{2\pi \textup{i} \eps\omega}|=1.
\end{equation}
It is shown in \cite{lawlor2013adaptive} that for a randomly chosen $\eps>0$ the
converse holds with probability one, and furthermore checking the condition
~\eqref{eq:non-collision}
for several $\eps$ would be sufficient to deterministically
decide whether $\oo$ has a collision. In section~\ref{sec:multi} we use this latter observation to devise a robust test for collisions even in the presence of noise. 

The algorithm developed in \cite{lawlor2013adaptive} for recovering $S(t)$ is as follows:
First we pick a prime $p=p_1$, which is roughly $5k$ where $k= |\OO|$ is the number of modes
in $S(t)$ ($k$ is commonly referred to as the {\em sparsity} of $S(t)$). By taking $p \ge 5k$ we ensure that on average collisions do not occur
for more than 90\% of $\oo \in\OO$. Let $\OO'$ denote the subset of $\OO$ consisting of all
non-collision $\oo\in\OO$. For each $\oo\in\OO'$ we recover $a_\oo \e^{2\pi \textup{i}\omega t}$, and
update $S(t)$ to
\begin{equation}\label{eq:s1}
   S_1(t) = S(t) - \sum_{\oo\in\OO'}a_\oo \e^{2\pi \textup{i}\omega t}.
\end{equation}
We now apply the above procedure again for $S_1(t)$ with a different prime $p=p_2$
approximately in the range of $5k_1$, where $k_1= k-|\OO'|$ is now the sparsity for
$S_1(t)$. This process is repeated until all modes are found.

In the implementation of the algorithm we set a small threshold in~\eqref{eq:non-collision}
to check for collisions. This means there is a small probability that a
collision is undetected by our criterion and a false value $\oo_0$ is put into $\OO'$ when it shouldn't be. In subsequent iterations, this will create a new mode $-c_0 \e^{2\pi \textup{i}\omega_0 t}$ for some $c_0\in\C$ in $S_1(t)$. By the use of different primes
$p_j$ in each iteration this false mode will very likely be identified and subtracted from the final reconstruction. In Subsection~\ref{subsec:multi_algo} we provide an improved aliasing test for our multiscale algorithm which makes the inclusion of spurious frequencies even less likely. However, it is still possible that incorrect modes are inserted and deleted in the high-noise regime, as we discuss in Section~\ref{sec:empirical}. 

\subsection{Noise model}
\label{subsec:noise_model}

In a number of potential application areas for sparse Fourier algorithms, the samples collected will be corrupted by noise. One example of sparse Fourier transforms being used on real data is given in \cite{hassanieh2012faster}, where an application to faster GPS location is presented. Previous works to address the issue of noise in the sparse Fourier transform context include \cite{hassanieh2012nearly}, although the algorithm presented in that work for noisy signals has yet to be implemented and evaluated empirically.

In this paper we assume an i.i.d. noise model
\begin{equation}  \label{eq:noisy_meas}
     \bm{S}^\mathrm{n}_p[j] = S\left(\frac{j}{p}\right)
       + \rmn_j = \bS_p[j] + \rmn_j,
\end{equation}
where $\rmn=(\rmn_j)$ are i.i.d. complex random variables with mean 0 and variance $\sigma^2$.
A typical model is to assume $\{\rmn_j\}$ are i.i.d. complex Gaussian. With the noise model we have
\begin{equation}  \label{noisy_dft}
   \wh{\bS}^\mathrm{n}_p[h] = \wh\bS_p[h] + \wh{\rmn}[h],
\end{equation}
where
\begin{equation}\label{eq:nhat}
   \wh{\rmn}[h] = \sum_{j=0}^{p-1} \rmn_j\e^{-2\pi\textup{i} hj/p}.
\end{equation}

By the i.i.d. property for $\{\rmn_j\}$ we have for each $h$
\begin{equation}\label{eq:en}
 \E\bigl[\wh{\rmn}[h]\bigr] = 0
\end{equation}
and
\begin{equation}\label{eq:varn}
   \Var\bigl[\wh{\rmn}[h]\bigr] = p \sigma^2,
\end{equation}
where the expectations are taken with respect to the randomness in the noise. This yields
\begin{equation}  \label{eq:mean}
    \E\left[\wh{\bS}^\mathrm{n}_p[h]\right] = \wh{\bS}_p[h]
\end{equation}
and
\begin{equation} \label{eq:mean_variance}
   \E\left[|\wh{\bS}^\mathrm{n}_p[h] - \wh{\bS}_p[h]|^2\right] = p \sigma^2.
\end{equation}
Thus, a typical noisy DFT coefficient $\wh{\bS}^\mathrm{n}_p[h]$ will deviate from the true
value $\wh{\bS}_p[h]$ by an amount proportional to $\sigma\sqrt{p}$. Similarly,
for $\bS_{p,\eps}^\rmn = \bS_{p,\eps}+ \rmn_\eps$ we will have
\begin{equation}  \label{eq:noisy_mean}
    \E\left[\wh{\bS}^\mathrm{n}_{p,\eps}[h]\right] = \wh{\bS}_{p,\eps}[h]
\end{equation}
and
\begin{equation} \label{eq:noisy_variance}
   \E\left[|\wh{\bS}^\mathrm{n}_{p,\eps}[h] - \wh{\bS}_{p,\eps}[h]|^2\right] = p \sigma^2.
\end{equation}

We now pick a non-collision $\oo \in\OO$. Then for $h = \oo \wmod p$ we will have
\begin{align} %\label{eq:pnoise_dft}
   \wh{\bS}^\mathrm{n}_p[h] &= p a_{\oo}  + \bigo(\sqrt p \sigma), \notag \\
   \wh{\bS}^\mathrm{n}_{p,\eps}[h] &= p a_{\oo}\e^{2\pi \textup{i} \oo\eps}  + \bigo(\sqrt p \sigma).
\end{align}
As a result $a_{\oo}$ can now be estimated easily via
\begin{equation} \label{eq:pnoise_aj}
    a_{\oo} = \frac{1}{p}\wh{\bS}^\mathrm{n}_p[h]  + \bigo\Bigl(\frac{\sigma}{\sqrt p}\Bigr).
\end{equation}

The real challenge lies in the recovery of the frequencies in $\OO$. Assume that
$|\wh{\bS}_{p,\eps}|$ has a pulse at $h$. Then $h = \oo \wmod p$ for some $\oo\in\OO$.
If there is no collision for $\oo$, in the noiseless environment $\oo$ is recovered
via~\eqref{eq:omega} as long as $\eps \le \frac{1}{N}$. In the noisy setting
$\wh{\bS}_{p,\eps}[h]/\wh{\bS}_p[h]$ must be replaced by
$\wh{\bS}^\mathrm{n}_{p,\eps}[h]/\wh{\bS}^\mathrm{n}_p[h]$. Interestingly,
the mean of $\wh{\bS}^\mathrm{n}_{p,\eps}[h]/\wh{\bS}^\mathrm{n}_p[h]$ is in general
{\em not} $\wh{\bS}_{p,\eps}[h]/\wh{\bS}_p[h]$ as a result of the division. Nevertheless
we have
\begin{align}
\label{eq:noisy_binom}
\frac{\wh{\bS}^\mathrm{n}_{p,\eps}[h]}{\wh{\bS}^\mathrm{n}_p[h]} &=
        \frac{\wh{\bS}_p[h]\e^{2\pi\textup{i} {\oo}\eps} + \wh\rmn_\eps[h]}
        {\wh{\bS}_p[h]+\wh\rmn[h]} \notag \\
      &= \frac{\wh{\bS}_p[h]\e^{2\pi\textup{i} {\oo}\eps} + \bigo\left(\sigma\sqrt{p}\right)}
        {\wh{\bS}_p[h]+\bigo\left(\sigma\sqrt{p}\right)} \notag \\
      &= \frac{\e^{2\pi\textup{i}\oo\eps} + \bigo\left(\sigma/a_{\oo}\sqrt{p} \right)}
      {1 + \bigo\left( \sigma/a_{\oo} \sqrt{p}\right)} \notag \\
      &= \e^{2\pi\textup{i}\oo \eps} + \bigo\left(\sigma/a_{\oo}\sqrt{p}\right).
\end{align}
Thus the ratio of noisy DFT coefficients agrees with the noiseless ratio up to an error term
on the order of $\sigma/|a_{\oo}|\sqrt{p}$.

Given this estimate for the ratio of noisy DFT coefficients, we can derive bounds for the
error in the Lee norm for the phase angle computed via $\Arg(z)$. Let ${\mathcal L}$
be a lattice in $\R$.
For any $\theta\in\R$ the {\em Lee norm associated with the lattice ${\mathcal L}$}
 for $\theta$ is given by the distance of $\theta$ to the lattice ${\mathcal L}$, i.e.
$\|\theta\|_{\mathcal L} := \min_{k\in {\mathcal L}} |\theta-k|$. Under the Lee norm
associated with the lattice $2\pi\Z$ it is well known that
\begin{align} \label{eq:arg_bound}
    \|\Arg\left(z+\eta\right) - \Arg(z)\|_{2\pi \Z} &=
     \|\Arg\left(1+z^{-1}\eta\right)\|_{2\pi \Z} \notag \\ &\le |z^{-1}\eta|.
\end{align}
Thus for a non-collision $\oo\in\OO$ and $h = \oo \wmod p$, the estimates~\eqref{eq:arg_bound} and~\eqref{eq:noisy_binom} combined yield
\begin{equation} \label{eq:phase_error}
    \left\|\Arg\left(\frac{\wh{\bS}^\mathrm{n}_{p,\eps}[h]}
      {\wh{\bS}^\mathrm{n}_p[h]}\right) - 2\pi\omega\eps
    \right\|_{2\pi \Z} \le \bigo\left( \frac{\sigma}{|a_{\oo}|\sqrt{p}}\right).
\end{equation}
When we apply the estimate~\eqref{eq:omega} for $\omega$ under the noise model we will
end up with an approximation
\begin{equation}\label{eq:omegan}
    \omega^\rmn := \frac{1}{2\pi  \eps}\Arg\left(\frac{\widehat{\bS}^\rmn_{p,\eps}[h]}
       {\widehat{\bS}^\rmn_{p}[h]}\right)
\end{equation}
such that
\begin{equation} \label{eq:noisy_omega}
     \|\omega^\rmn-\omega\|_{\Z} \le \bigo\left(\frac{\sigma}{2\pi\eps|a_{\oo}|\sqrt{p}}\right).
\end{equation}
Now if we apply the algorithm developed in \cite{lawlor2013adaptive} the ratio
$\frac{\sigma}{\eps\sqrt{p}}$ is critical in determining the sensitivity of our phase estimation (as
well as the weight estimation) to noise. Without any modifications to the algorithm it is thus
important that we choose the lengths $p$ so that $\frac{\sigma}{\eps\sqrt{p}}$ is within the
tolerance.

\subsection{Earth mover distance}
\label{subsec:emd}

In the existing literature on the sparse Fourier transform, the $\ell_2$ norm is most often used to assess the quality of approximation. There are many reasons for this choice, with the two most convincing perhaps being the completeness of the complex exponentials with respect to the $\ell_2$ norm and Parseval's theorem. For certain applications, however, this choice of norm is inappropriate. For example, in wide-band spectral estimation and radar applications, one is interested in identifying a set of frequency intervals containing active Fourier modes. In this case, an estimate $\wt{\omega}$ of the true frequency $\omega$ with $|\wt{\omega}-\omega| \ll N$ is useful, but unless $\wt{\omega}=\omega$ the $\ell_2$ metric will report an $\bigo(1)$ error.
For these reasons, we propose measuring the approximation error of sparse Fourier transform problems with the Earth Mover Distance (EMD) \cite{rubner2000earth}. Originally developed in the context of content-based image retrieval, EMD measures the minimum cost that must be paid (with a user-specified cost function) to transform one distribution of points into another. EMD can be calculated efficiently as the solution of a linear program corresponding to a certain flow minimization problem. 

For our problem, we consider the cost to move a set of estimated Fourier modes and coefficients $\left\{(\wt{\omega}_j,a_{\wt{\omega}_j})\right\}_{j=1}^{\wt{k}}$ to the true values $\left\{(\omega_j, a_{\omega_j})\right\}_{j=1}^k$ under the cost function
\begin{equation}\label{eq:d1}
d_1\big( (\omega,a_\omega), (\wt{\omega},a_{\wt{\omega}}); N \big) := \frac{|\omega-\wt{\omega}|}{N} + |a_\omega-a_{\wt{\omega}}|.
\end{equation}
This choice of cost function strikes a balance between the fidelity of the frequency estimate (as a fraction of the bandwidth) and that of the coefficient estimate. We also consider the ``phase-only'' cost function 
\begin{equation} \label{eq:emd-phase}
d_\omega(\omega_1,\omega_2; N) := \frac{|\omega_1-\omega_2|}{N},
\end{equation}
which provides a measure of how close our frequency estimates are to the true values. We denote the EMD using $d_1$ by EMD(1) and using $d_\omega$ by EMD($\omega$) in our empirical studies in Section~\ref{sec:empirical} below.

Since these error metrics may be unfamiliar to the reader, we note here that the theoretical best possible EMD(1) error is easy to compute in the special case when the EMD($\omega$) error is zero (i.e., all frequencies are estimated correctly). In this case, we can combine~\eqref{eq:pnoise_aj} with~\eqref{eq:d1} above to yield
\begin{equation}\label{eq:emd1-specialcase}
  \textup{EMD}(1) = \bigo\left(\frac{k\sigma}{\sqrt{p}} \right).
\end{equation}
Note in particular that since we measure distances in $\ell_1$ the error scales with $k$, rather than $\sqrt{k}$ as would be the case in $\ell_2$. The case when EMD($\omega$) is non-zero is much more difficult to analyze and is an important question that merits considerable attention. We plan to conduct such a study in future work.

\subsection{Random signal model}
\label{subsec:random_signal_model}
For the empirical evaluations in Section~\ref{sec:empirical} we consider test signals with uniformly random phase over the bandwidth and coefficients chosen uniformly from the complex unit circle. In other words, given $k$ and $N$, we choose $k$ frequencies $\omega_j$ uniformly at random (without replacement) from $[-N/2,N/2)\cap\Z$. The corresponding Fourier coefficients $a_j$ are of the form $\e^{2\pi\textup{i}\theta_j}$, where $\theta_j$ is drawn uniformly from $[0,1)$.  The signal is then given by
\begin{equation}
    \label{eq:test_signal}
    S(t) = \sum_{j=1}^k a_j \e^{2\pi\textup{i}\omega_j t}.
\end{equation}
This is the standard signal model considered in previous empirical evaluations of sub-linear Fourier algorithms \cite{iwen2007empirical, iwen2010combinatorial, hassanieh2012simple, lawlor2013adaptive}. We note here that we also conducted the empirical evaluations of Section~\ref{sec:empirical} on signals whose Fourier coefficients have varying magnitudes. These results did not differ substantively from those on signals of the form~\eqref{eq:test_signal}, so we omit a detailed discussion.

%%%%%%%%%%%%%%%%%%%%%%%%%%%%%%
\section{Rounding: A Minor Modification of Noiseless Algorithm}
\label{sec:minor}

A simple modification to the noiseless algorithm of \cite{lawlor2013adaptive} for the noisy
case is to increase the sample lengths $p$. By choosing $p$ large enough the error from noise
can be mitigated to be within a given tolerance. The modification can be viewed simply as
rounding, and we include it both as a more direct and simple to implement extension as well as
for comparison purposes. When the noise level is low, this modification yields reasonably good
results. 

As in the noiseless case we choose the shift $\eps>0$ so that $\eps \leq \frac{1}{N}$. In the
noiseless case $\eps = \frac{1}{N}$ would be sufficient to avoid wrap-around aliasing in the phase
reconstruction. Due to the presence of noise we will need to make $\eps$ slightly smaller
because of~\eqref{eq:noisy_omega}. Let us analyze the recovery of a candidate frequency $\oo\in\OO$ if
we simply carry out the same process as in the noiseless environment.

First we choose a length $p$. Assume that $\oo\in\OO$ does not collide with any other
$\oo'\in\OO$ modulo $p$. Let $h = \omega \wmod p$. The reconstruction of $\omega$ utilizes two
factors. First, the location of peaks in the DFT are robust to noise: even with a relatively high
noise level we may take $h=\omega \wmod p$ to be exact. Second, by~\eqref{eq:noisy_omega} the
frequency reconstruction from noisy measurements is correct up to an error term of size
$\bigo\left( \frac{\sigma}{2\pi\eps|a_\oo|\sqrt{p}}\right)$. By combining these two measures we can
more reliably estimate $\omega$.

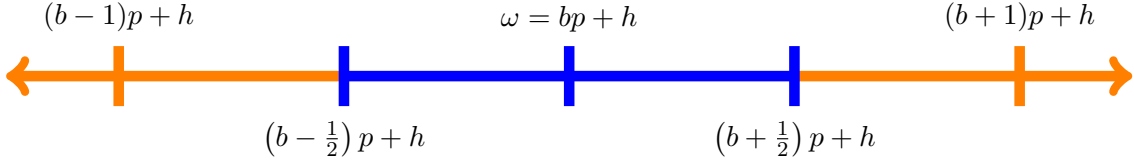
\begin{figure}
  \centering
  \resizebox{\columnwidth}{!}{%
    \begin{tikzpicture}[xscale=15]
      \draw[<-][draw=orange,line width = 4] (0,0) -- (0.3,0);
      \draw[draw=blue,line width = 4](0.3,0) -- (0.7,0);
      \draw[->][draw=orange,line width = 4](0.7,0) -- (1,0);
      \draw[draw=blue,line width = 4] (0.5,-0.4) -- (0.5,0.4) node[above]{$\omega = bp+h$}; 
      \draw[draw=orange,line width = 4] (0.1,-0.4)  -- (0.1,0.4) node[above]{$(b-1)p+h$}; 
      \draw[draw=orange,line width = 4] (0.9,-0.4)  -- (0.9,0.4) node[above]{$(b+1)p+h$};
      \draw[draw=blue,line width=4] (0.3,-.4) node[below]{$\left(b-\tfrac{1}{2}\right)p+h$} -- (0.3,0.4);
      \draw[draw=blue,line width=4] (0.7,-.4) node[below]{$\left(b+\tfrac{1}{2}\right)p+h$} -- (0.7,0.4);
    \end{tikzpicture}
  }
  \caption{The rounding procedure is exact as long as the phase estimate $\wt{\omega}$ is within $p/2$ of correct multiple of $p$ (blue region in figure).}
  \label{fig:rounding}
\end{figure}

Our proposed modification is to simply round the noisy frequency estimate
\begin{equation}\label{eq:omegatilde}
    \wt{\omega}= \frac{1}{2\pi\eps}\Arg\left( \frac{\wh{\bS}^\mathrm{n}_{p,\eps}[h]}{\wh{\bS}^\mathrm{n}_p[h]} \right)
\end{equation}
to the nearest integer of the form $np+h$. This improved estimate is therefore given by
\begin{equation}
\label{eq:round_omega}
\wt{\omega}' = p\cdot \mathrm{round}\left( \frac{\wt{\omega}-h}{p}\right) + h,
\end{equation}
where $\mathrm{round}(x)$ returns the nearest integer to $x$. For low noise levels this
modification will return the true value $\omega$, while for larger noise levels it is possible
that $\wt{\omega}$ deviates by more than $p/2$ from the true frequency $\omega$. In this case
the estimate $\wt{\omega}'$ will be wrong by a multiple of $p$. Larger values of $p$ will
reduce the likelihood of an error in frequency estimation. See Figure~\ref{fig:rounding} for an illustration of this rounding procedure.

To ensure that the estimated frequencies are sufficiently far from the branch cut of
$\Arg(z)$ along the negative real axis, we take the shift $\eps \le \frac{1}{2N}$. The estimated
frequencies then satisfy $-N \le \wt{\omega} < N$, while the true frequencies lie in the smaller interval $[-N/2,N/2)$. It is thus extremely unlikely that the deviations due to the noise
will push the estimates across the discontinuity.

We saw in the previous section that the error in the phase estimation is on the order of
$\sigma p^{-1/2}$ when using the reconstruction formula~\eqref{eq:omega}. When using the
rounding procedure~\eqref{eq:round_omega}, however, we should expect accurate results for a
wider range of sample lengths $p$ and noise levels $\sigma$. Indeed, note that the rounded
frequency estimate $\wt{\omega}'$ is \emph{exact} as long as
\begin{equation}
\label{eq:exact_round}
|\wt{\omega}-\omega| < \frac{p}{2}.
\end{equation}
Recall from Section~\ref{subsec:noise_model} that the error of the frequency estimate
$\wt{\omega}$ is on the order of $\bigo(\frac{\sigma}{\eps\sqrt{p}})$. Let us assume that it is
bounded by $C\frac{\sigma}{\eps\sqrt{p}}$ for some constant $C$. Combining this with the requirement~\eqref{eq:exact_round} we see that the rounded frequency estimate $\wt{\omega}'$ will be exact
provided
\begin{equation}
\label{eq:round_req}
    C\frac{\sigma}{\eps p^{3/2}} < \frac{1}{2}.
\end{equation}
It follows that we get exact reconstruction if $p \geq (2C\sigma/\eps)^{2/3}$.

\begin{figure}
    \centering
    \includegraphics[width=0.9\columnwidth]{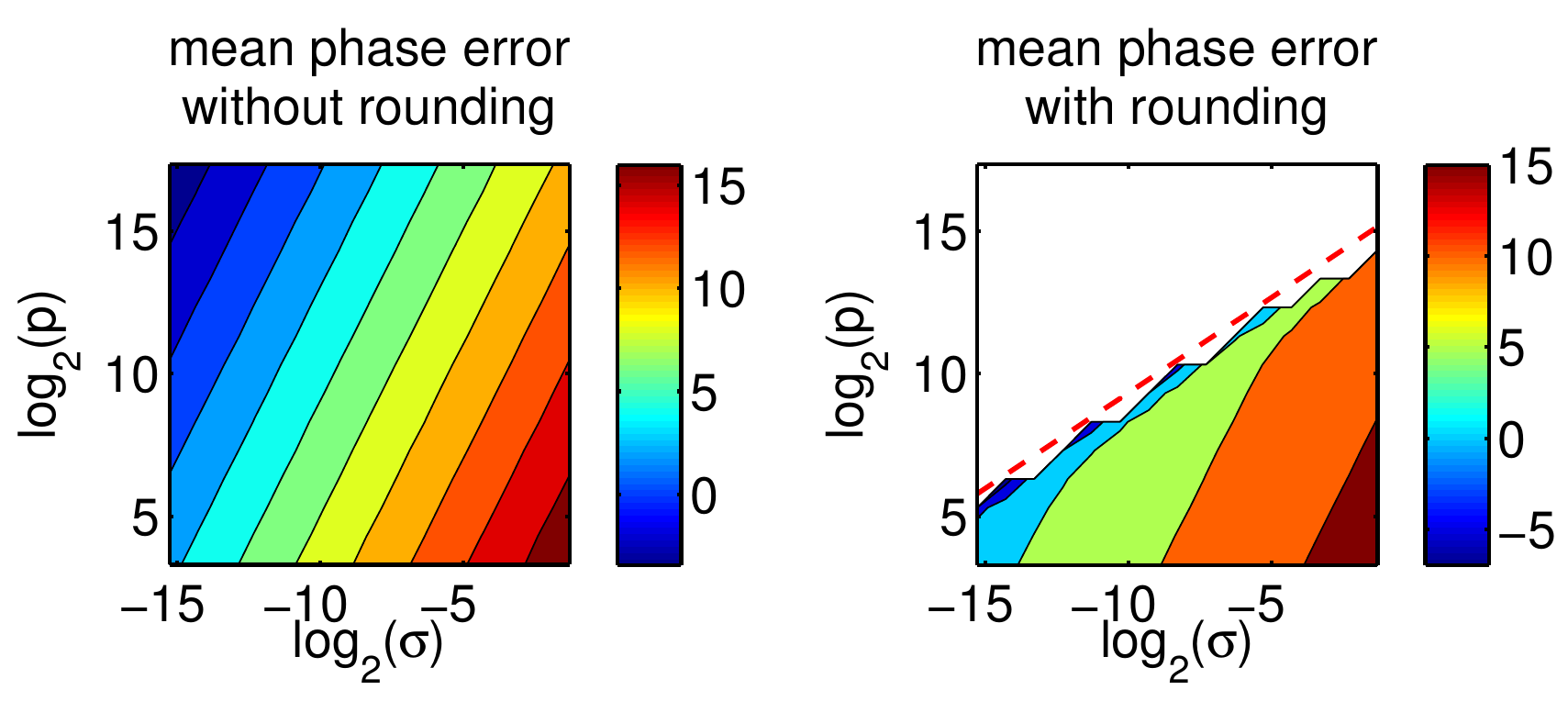}
    \caption{(left) Mean phase error (in log scale) for frequency estimation via~\eqref{eq:omega}. (right) Mean phase error (in log scale) for frequency estimation with rounding via~\eqref{eq:round_omega}. The red dotted line marks the transition to exact recovery when $p > (2\sigma/\eps)^{2/3}$.}
    \label{fig:phase_error}
\end{figure}

To illustrate this relationship, we generated 1000 test signals with
frequencies chosen uniformly at random from $[-N/2, N/2)$ and set the
  corresponding coefficient to unity. Thus our test signals for this
  empirical trial were one-term trigonometric polynomials. For this
  test we took $N=2^{22}$, $\eps=\frac{1}{2N}$ and investigated a range of parameters
  $(\sigma, p)$.  We reconstructed the frequencies in two ways: first,
  simply using the formula~\eqref{eq:omega}, and second by combining
  this estimate with the rounding procedure~\eqref{eq:round_omega}. In
  Figure \ref{fig:phase_error} we plot the average phase error in
  logarithmic scale as a function of both $\sigma$ and $p$, which were
  varied from $2.5\times10^{-5}$ to $0.4096$ and from $10$ to
  $163840$, respectively, by powers of two.

In the plot on the left, which corresponds to reconstruction using only~\eqref{eq:omega}, we can clearly see the contours of constant phase error obeying the relationship $\log_2(p) = 2 \log_2(\sigma) + \alpha$ for various $\alpha$. This confirms our analytic estimate from Subsection~\ref{subsec:noise_model} that the phase error is proportional to $\sigma/\sqrt{p}$. In the plot on the right, which corresponds to the improved reconstruction using~\eqref{eq:round_omega}, we can see that for large values of $\sigma$ and small values of $p$ the same relationship holds. However, for smaller $\sigma$ and larger $p$ we see an abrupt transition to exact reconstruction (the white area in the upper-left). The boundary of this region (red dotted line) follows the relationship $\log_2(p) = \frac{2}{3}\log_2(\sigma) + 16$, corresponding to $C=1$ in~\eqref{eq:round_req} above. This illustrates that for small enough values of the ratio $\frac{\sigma}{\eps p^{3/2}}$ the rounding procedure is exact.

\subsection{Algorithm}
\label{subsec:round_algo}
Our first algorithm for noisy signals is only a slight modification of the noiseless algorithm presented in \cite[Algorithm 1]{lawlor2013adaptive}. Considering~\eqref{eq:round_req}, we change the lower bound
\begin{equation}
\label{eq:old_p}
p>c_1 k
\end{equation}
to
\begin{equation}
\label{eq:new_p}
p>\max\{c_1k, c_2(\sigma/\eps)^{2/3}\},
\end{equation}
where $c_1,\, c_2$ are constants.
In this way we ensure that the choice of $p$ is always large enough to isolate most of the $k$ frequencies on average as well as being large enough to ensure that the rounding procedure~\eqref{eq:round_omega} is exact. In all of our experiments in Section~\ref{sec:empirical} below we took $c_2=4$.

%%%%%%%%%%%%%%%%%%%%%%%%%%%%%%%%%%%%%%%%%%%%%%
\section{A Multiscale Algorithm}
\label{sec:multi}

In Section~\ref{sec:minor} we saw that taking $p > \max\{c_1k,c_2(\sigma/\eps)^{2/3}\}$ sufficed to ensure that the rounding procedure was exact. While this gives good results in terms of accuracy, the increased runtime associated with larger noise levels is undesirable. The main contribution of this paper is a multiscale algorithm for recovering the
frequency set $\Omega$ of the signal $S(t)$. This algorithm achieves similar accuracy while providing improvement by several orders of magnitude in computational efficiency.

The key feature of this multiscale algorithm is the employment of
multiple shifts $\eps_j$, which
enable us to improve the accuracy of the phase estimations progressively without the need
to significantly increase the sample length $p$. As we will see, taking successively larger shifts enables a form of error-correction in our frequency estimates at finer and finer scales, in essence ``zooming in'' on the true frequencies in a multiscale fashion. In Subsection~\ref{subsec:multi_freqid} we give some background on our multi-scale method and introduce the main idea of our algorithm. In Subsection~\ref{subsec:multi_tech} we prove that our multiscale approximations are accurate estimates of the true frequencies, and in Subsection~\ref{subsec:multi_algo} we describe the basic multiscale algorithm.

\subsection{Multiscale frequency estimation}
\label{subsec:multi_freqid}
The main idea for the multiscale algorithm is that a value can be estimated with high precision with an inaccurate (coarse) estimator applied progressively at different
scales, much like in analog-to-digital
conversion where a signal value can be estimated with very high precision by the very coarse
binary quantization. In our sparse Fourier recovery algorithm, the coarse estimator is the
approximation formula given by~\eqref{eq:phase_error}
\begin{equation}
    \eps\omega =_{\Z} \frac{1}{2\pi}\Arg\left(\frac{\wh{\bS}^\mathrm{n}_{p,\eps}[h]}
      {\wh{\bS}^\mathrm{n}_p[h]}\right),
\end{equation}
where $=_{\Z}$ is measured by the Lee norm $\|\cdot\|_\Z$.

For simplicity let us assume for the moment that our signal contains a single frequency $\omega$ with non-zero Fourier coefficient. For a fixed $p$, let $\wt{\omega}$ be our estimate for $\omega$ using the rounding procedure from Section~\ref{sec:minor} with shift $\eps_0 \le \frac{1}{N}$. Then we have
\begin{equation}
\label{eq:omega_omegatilde}
\wt{\omega} = \omega \wmod{p},
\end{equation}
although in general $\wt{\omega}$ may differ from $\omega$ by a multiple of $p$.

Suppose now that we repeat the computation of $\wt{\omega}$ using a larger shift $\eps_1 > \eps_0$; that is, we sample our signal at time points $\frac{j}{p} + \eps_1$, take the FFT, and compute
\begin{equation}
\label{eq:b_1}
b_1 = \frac{1}{2\pi}\Arg\left( \frac{\wh{\bS}^\mathrm{n}_{p,\eps_1}[h]}{\wh{\bS}^\mathrm{n}_p[h]} \right)
\end{equation}
(note that we do not divide by $\eps_1$). Since in general $\eps_1 > \frac{1}{N}$, we cannot take $b_1/\eps_1$ as an estimate for $\omega$, although it still holds that
\begin{equation}
\label{eq:b_1mod1}
b_1 \approx \eps_1 \omega \wmodc,
\end{equation}
where $x \wmodc$ is the unique value $y$ in $[-\frac{1}{2},\frac{1}{2})$ such that $x\equiv y \wmod 1$. We can use this fact to estimate the error $\omega-\wt{\omega}$ as follows. Note that
\begin{align}
\label{eq:eps_diff_omega}
\eps_1(\omega-\wt{\omega}) & = \eps_1\omega -\eps_1\wt{\omega} \\ \notag
& \approx (b_1 - \eps_1\wt{\omega}) \wmodc,
\end{align}
so that
\begin{equation}
\label{eq:diff_omega}
\omega -\wt{\omega} \approx (b_1 - \eps_1\wt{\omega})\wmodc /\eps_1.
\end{equation}

\begin{figure}
  \centering
  \resizebox{\columnwidth}{!}{%
    \begin{tikzpicture}[xscale=15]
      
      \draw[draw=blue,line width = 10] (0,0) -- (0.25,0);
      \draw[draw=orange,line width = 10](0.25,0) -- (0.5,0);
      \draw[draw=orange,line width = 10](0.5,0) -- (0.75,0);
      \draw[draw=orange,line width = 10](0.75,0) -- (1,0) node[right]{iteration 1};

      \draw[draw=blue,line width = 2] (0,-0.4) -- (0,0.4) node[above]{$\eps_0^{-1}$}; 
      \draw[draw=orange,line width = 2] (0.25,-0.4) -- (0.25,0.4) node[above]{$\eps_1^{-1}$}; 
      \draw[draw=orange,line width = 2] (0.5,-0.4) -- (0.5,0.4) node[above]{$\cdots$};
      \draw[draw=orange,line width = 2] (0.75,-.4) -- (0.75,0.4) node[above]{$\eps_m^{-1}$};
      \draw[draw=orange,line width = 2] (1,-.4) -- (1,0.4) node[above]{$1$};

      \draw[draw=blue,line width = 10](0.25,-1) -- (0.5,-1);
      \draw[draw=orange,line width = 10](0.5,-1) -- (0.75,-1);
      \draw[draw=orange,line width = 10](0.75,-1) -- (1,-1) node[right]{iteration 2};
      \draw[draw=blue,line width = 2] (0.25,-1.4) -- (0.25,-0.6);
      \draw[draw=orange,line width = 2] (0.5,-1.4) -- (0.5,-0.6);
      \draw[draw=orange,line width = 2] (0.75,-1.4) -- (0.75,-0.6);
      \draw[draw=orange,line width = 2] (1,-1.4) -- (1,-0.6);

      \draw[draw=blue,line width = 10](0.75,-3) -- (1,-3) node[right]{iteration $m$};
      \draw[draw=blue,line width = 2] (0.75,-3.4) -- (0.75,-2.6) node[above]{$\vdots$};
      \draw[draw=blue,line width = 2] (1,-3.4) -- (1,-2.6) node[above]{$\vdots$};

    \end{tikzpicture}
  }
  \label{fig:multi-diagram}
  \caption{Diagram of the multiscale frequency estimation procedure, with a candidate frequency pictured as a string of digits, from most significant on the left to least significant on the right. In this figure, blue regions represent correct digits learned by the algorithm, and orange regions represent digits where errors are likely. In the first iteration, the most significant bits are learned using shift $\eps_0^{-1}$. Subsequent iterations give corrections at finer scales $\eps_1^{-1}, \ldots, \eps_m^{-1}$.}
\end{figure}
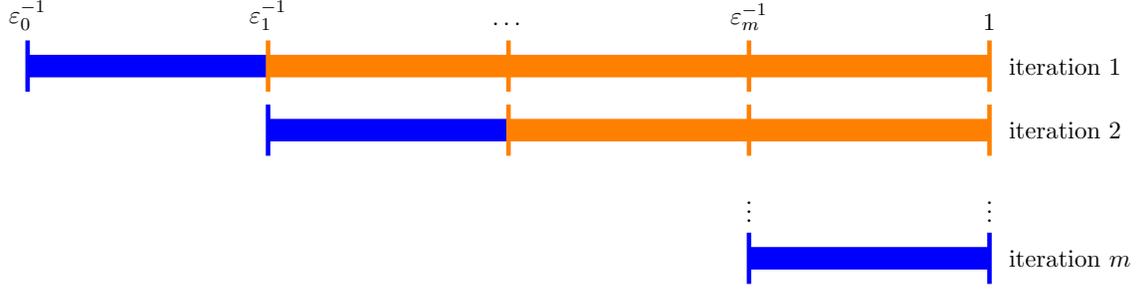

This estimate of the error is not exact, since there is still noise that can perturb the calculated value $b_1$ from the true value $\eps_1 \omega \wmodc$. However, analogous to~\eqref{eq:noisy_omega} we have
\begin{equation}\label{eq:omega_error}
(\omega-\wt{\omega}) - (b_1 - \eps_1\wt{\omega})\wmodc/\eps_1 = \bigo\left(\frac{\sigma}{\eps_1\sqrt{p}}\right),
\end{equation}
which immediately implies that the updated estimate satisfies
\begin{equation}
\label{eq:update_omega_error}
\omega - \left(\wt{\omega}+(b_1-\eps_1\wt{\omega})\wmodc/\eps_1\right) = \bigo\left(\frac{\sigma}{\eps_1\sqrt{p}}\right).
\end{equation}
Since $\eps_1 > \eps_0$, adding the correction term~\eqref{eq:diff_omega} to our previous estimate $\wt{\omega}$ will give a finer approximation to the true frequency $\omega$. By iterating this error correction process with progressively larger shifts $\eps_j$, we obtain an algorithm which adaptively corrects for the error in a multiscale fashion. See  Figure~\ref{fig:multi-diagram} for a diagram of the multiscale estimation procedure. In the next section we provide a detailed analysis of this multiscale approximation scheme, and prove that the frequency estimates it produces are accurate.

\subsection{Analysis of multiscale approximations}
\label{subsec:multi_tech}
We begin with a technical lemma relating arithmetic in the Lee norm $\|\cdot\|_\Z$ to that on the interval $[-\frac{1}{2},\frac{1}{2})$. It will be used repeatedly in the sequel.

\begin{lem}  \label{lem-4.1}
   Let $\delta>0$ and $x\in [-\frac{1}{2}+\delta, \frac{1}{2}-\delta]$. Assume that $\|x-b\|_\Z <\delta$
   and $b\in [-\frac{1}{2}, \frac{1}{2})$. Then $|x-b| <\delta$.
\end{lem}
\Proof  Let $r=\|x-b\|_\Z$. Then $x-b =\pm r+k$ for some $k\in\Z$. If $k=0$ we have
\begin{equation}\label{eq:xmb}
 |x-b| = \| x-b \|_{\Z} < \delta
\end{equation}
by hypothesis, so the claim holds. Now assume $k\neq 0$. Note that 
\begin{equation}
\label{eq:case1}
|x-b| \leq |x|+|b| \leq 1-\delta
\end{equation}
by the triangle inequality and the assumptions on $x$ and $b$. At the same time, we have
\begin{equation}
\label{eq:case2}
|\pm r+k|\geq 1-r >1-\delta.
\end{equation}
This is a contradiction, since~\eqref{eq:case1} and~\eqref{eq:case2} cannot hold simultaneously. Thus we must in fact have $k=0$, and the claim holds.
\eproof

The following Theorem formalizes the multiscale frequency estimation procedure which was introduced in the previous Subsection.
\begin{theo}  \label{theo-4.2}
   Let $\oo\in [-\frac{N}{2}, \frac{N}{2})$. Let $0 <\eps_0 <\eps_1 < \cdots <\eps_m$ and
   $b_0, b_1, \dots, b_m\in\R$ such that
   \begin{equation}\label{eq:thm4.2hyp}
       \|\eps_j \oo - b_j\|_\Z < \delta, \mhsp 0 \leq j \leq m
   \end{equation}
   where $0<\delta \leq \frac{1}{4}$. Assume that $\eps_0 \leq \frac{1-2\delta}{N}$
   and $\beta_j:=\eps_j/\eps_{j-1} \leq (1-2\delta)/(2\delta)$. Then there exist
   $c_0, c_1, \dots, c_m\in\R$, each computable from $\{\eps_j\}$ and $\{b_j\}$,
   such that
   \begin{equation}\label{eq:thm4.2conc}
        |\wt\oo -\oo| \leq \frac{\delta}{\eps_0} \prod_{j=1}^m\beta_j^{-1}, \mhsp \mbox{where}\mhsp
        \wt\oo := \sum_{j=0}^m \frac{c_j}{\eps_j}.
  \end{equation}
\end{theo}
\Proof
   Denote $\oo_0:=\oo$. We first note that  
\begin{equation}
|\eps_0\oo_0| \leq \eps_0 \frac{N}{2}\leq \frac{1}{2}-\delta,
\end{equation}
where the second inequality follows from the assumptions of the Theorem. Let $c_0 = b_0 \wmodc$, so that $|\eps_0\oo_0-c_0|<\delta$ by Lemma~\ref{lem-4.1}. Let $\lambda_0 = c_0/\eps_0$, which represents a coarse estimate of $\oo_0$ with the error bound 
\begin{equation}
|\lambda_0-\oo_0| < \delta/\eps_0.
\end{equation}

Next, let $\oo_1=\oo_0-\lambda_0$. By the above  $|\oo_1| < \delta/\eps_0$ and
\begin{equation}
      |\eps_1\oo_1 |< \frac{\eps_1\delta}{\eps_0} =\beta_1\delta \leq\frac{1}{2} -\delta.
\end{equation}
We then have
\begin{equation}
\|\eps_1\oo - b_1\|_\Z=\|\eps_1\oo_1 - (b_1-\eps_1\lambda_0)\|_\Z <\delta.
\end{equation}
Set
$c_1 = b_1-\eps_1\lambda_0 \wmodc$. It follows from Lemma~\ref{lem-4.1} again that
$|\eps_1\oo_1 -c_1| < \delta$. We set $\lambda_1= c_1/\eps_1$.

We can recursively define $c_j, \lambda_j$ and $\oo_j$ for all $1 \leq j \leq m$.
In general we define $\oo_{j} := \oo_{j-1}-\lambda_{j-1}$. This leads to
\begin{equation}
      |\eps_{j}\oo_j|< \frac{\eps_{j}\delta}{\eps_{j-1}} =\beta_j\delta \leq\frac{1}{2} -\delta.
\end{equation}
Set 
\begin{equation}
c_j = (b_j-\eps_j\lambda_{j-1}) \wmodc,
\end{equation}
which yields 
\begin{equation}
\|\eps_j\oo_j -c_j\|_\Z <\delta.
\end{equation}
Lemma~\ref{lem-4.1} now gives
$|\eps_j\oo_j - c_j| <\delta$. Set $\lambda_j = c_j/\eps_j$.

Finally denote $\oo_{m+1} = \oo_m-\lambda_m$. It is straightforward now to verify that
\begin{align}
      \oo=\oo_0 &= \sum_{j=0}^m \lambda_j +\oo_{m+1} \notag \\ &= \sum_{j=0}^m \frac{c_j}{\eps_j} +\oo_{m+1}.
\end{align}
Furthermore, by construction $\oo_{m+1} = \oo_m-\lambda_m$, which
has $|\oo_{m+1}| \leq \delta/\eps_m$. By hypothesis
$\eps_m = \eps_0\prod_{j=1}^m \beta_j$, yielding
\begin{equation}
     |\oo_{m+1}| \leq \frac{\delta}{\eps_0} \prod_{j=1}^m\beta_j^{-1}
\end{equation}
and completing the proof.
\eproof

\vspace{2mm}
\noindent
{\bf Remark 4.1.} ~From the proof of Theorem~\ref{theo-4.2} the values
$c_j$ and $\wt\oo$ are explicitly computable through the recursive formula $\oo_0=\oo$, $c_0 =
b_0 \wmodc$, $\lambda_0 = c_0/\eps_0$ and
\begin{equation} \label{eq:recursive}
       \left\{\begin{array}{ccl}
          \oo_{j} &=& \oo_{j-1}-\lambda_{j-1} \\
            c_j &=& (b_j-\eps_j\lambda_{j-1}) \wmodc \\
         \lambda_j &=& c_j/\eps_j
         \end{array}\right.
\end{equation}
for $1\leq j\leq m$. Equivalently, we can write the updated frequency estimates along the lines of~\eqref{eq:update_omega_error} as 
\begin{align} \label{eq:recursive-omega}
  \wt{\omega}_0 &= b_0/\eps_0 \notag \\
  \wt{\omega}_{n+1} &= \wt{\omega}_n + \left(b_n-\eps_n\wt{\omega}_n\right)\wmodc/\eps_n.
\end{align}

\vspace{2mm}
\begin{coro}   \label{coro-4.3}
      Assume that in the above theorem we have $\beta_j =\beta$ where
$\beta \leq (1-2\delta)/(2\delta)$, i.e. $\eps_j = \beta^j\eps_0$ for all $j$. Let $p>0$ and
$m \geq \left\lfloor\log_\beta \frac{2\delta}{p\eps_0}\right\rfloor+1$. Then
\begin{equation}
        |\wt\oo -\oo| \leq \frac{\delta}{\eps_0} \beta^{-m}< \frac{p}{2}.
\end{equation}
\end{coro}
\Proof   This is a straightforward corollary. By Theorem~\ref{theo-4.2} we have
\begin{equation}
        |\wt\oo -\oo| \leq \frac{\delta}{\eps_0} \prod_{j=1}^m\beta_j^{-1} = \frac{\delta}{\eps_0}\beta^{-m}.
\end{equation}
It is easy to check that $m =\left\lfloor\log_\beta \frac{2\delta}{p\eps_0}\right\rfloor+1$
is the smallest
integer such that $\frac{\delta}{\eps_0}\beta^{-m}<\frac{p}{2}$.
\eproof

Note that as we have mentioned in Section 3, even with noise the value $\oo \wmod p$ can be
precisely computed very reliably. Thus if the difference $|\oo-\wt\oo|$ is
smaller than $\frac{p}{2}$ then $\oo$ can be recovered exactly by taking
the closest integer to $\wt\oo$ with the same residue modulo $p$.

In numerical tests we choose uniform $\beta_j=\beta$. While making $\beta$ as large as
it can be for a given error estimate $\delta$ will undoubtedly reduce the computational cost,
there is nevertheless a good reason that we should not be too ``greedy'' and be more
conservative by choosing a smaller $\beta>1$. The reason is that given the random nature of
the noise the error bound $\delta$ is only in the average sense. To minimize reconstruction
errors we should try to provide as much latitude as possible for the uncertainties associated
with the error estimate $\delta$. Hence it is useful to ask how much latitude does one get for
given choices of $\eps_0$ and $\beta$.

\begin{theo}  \label{theo-4.4}
Let $\oo\in [-\frac{N}{2}, \frac{N}{2})$, $\eps_0>0$ and $\beta>1$. Set $\eps_j=\beta^j\eps_0$ for $1\leq j
\leq m$. Assume that we have $b_0, b_1, \dots, b_m\in\R$ such that
\begin{equation}
       \|\eps_j \oo - b_j\|_\Z < \delta, \mhsp 1 \leq j \leq m
\end{equation}
   where
   \begin{equation} \label{delta-bound}
       \delta = \min\,\Bigl(\frac{1-\eps_0N}{2}, \frac{1}{2\beta+2}\Bigr).
   \end{equation}
   Then the estimate $\wt\oo$ of $\oo$ given by $\wt\oo := \sum_{j=0}^m \frac{c_j}{\eps_j}$ satisfies
   \begin{equation}
        |\wt\oo -\oo| \leq \frac{\delta}{\eps_0} \beta^{-m},
   \end{equation}
   where $c_j$ are given in~\eqref{eq:recursive}.
\end{theo}
\Proof  The proof is straightforward. Note that Theorem~\ref{theo-4.2}
holds under the conditions $\eps_0 \leq \frac{1-2\delta}{N}$ and
$\beta_j \leq \frac{1-2\delta}{2\delta}$. These
two conditions are equivalent to the condition
$\delta \leq \min\,\bigl(\frac{1-\eps_0N}{2},
\frac{1}{2\beta+2}\bigr)$. Clearly, $\delta = \min\,\bigl(\frac{1-\eps_0N}{2},
\frac{1}{2\beta+2}\bigr)$ is the largest admissible value for $\delta$.
\eproof

\subsection{Algorithm}
\label{subsec:multi_algo}
 In this section we provide some details of our implementation of the multiscale frequency estimation procedure described in Subsection~\ref{subsec:multi_freqid}. In particular, we discuss the choice of various parameters necessary for reconstruction according to Theorem~\ref{theo-4.2} as well as changes made to the aliasing detection test from \cite{lawlor2013adaptive} to improve robustness in the presence of noise.

\subsubsection{Choice of $p$}
\label{subsubsec:choice-beta}
It remains to determine the choice of sampling length $p$, given the parameter $\beta$ and the noise level $\sigma$. Recall from the proof of Theorem~\ref{theo-4.2} that the estimated frequency $\wt\omega$ is given by the sum $\sum_{j=1}^m \lambda_j$, where $\lambda_j = c_j/\eps_j$. Moreover, the difference between successive frequency approximations is given in terms of $\lambda_j$ as
\begin{equation}
 \omega_j := \omega_{j-1}-\lambda_{j-1} \implies \lambda_j =\omega_j-\omega_{j+1}. 
\end{equation}
Thus we can decompose the error of approximation at stage $j+1$ as
\begin{align}
\label{eq:approx-error}
|\omega-\omega_{j+1}| &= |(\omega_j-\omega_{j+1})-(\omega_j-\omega)| \notag \\
&= |\lambda_j-(\omega_j-\omega)|.
\end{align} 

By Theorem~\ref{theo-4.2} the left-hand side of~\eqref{eq:approx-error} satisfies
\begin{equation}
 |\omega-\omega_{j+1}| < \frac{\delta}{\eps_{j+1}}, 
\end{equation}
while analogously to \eqref{eq:noisy_omega} the right-hand side of~\eqref{eq:approx-error} satisfies 
\begin{equation}
 |\lambda_j-(\omega_j-\omega)| \le \bigo\left(\frac{\sigma}{2\pi\eps_j\sqrt{p}}\right). 
\end{equation}
Denoting by $c_\sigma$ the constant in the right-hand side above and equating the two upper bounds gives 
\begin{equation}
\label{eq:choice-beta}
\frac{2\pi\delta\sqrt{p}}{c_\sigma \sigma} = \frac{\eps_{j+1}}{\eps_j} =: \beta.
\end{equation}

Under the assumptions of Theorem~\ref{theo-4.4}, we have
\begin{equation}
 \delta = \min\left(\frac{1-\eps_0 N}{2},\frac{1}{2\beta+2}\right).
\end{equation}
Since we take $\eps_0 = \frac{1}{2N}$ and fix $\beta > 1$, the latter term is necessarily the smaller. Plugging this into~\eqref{eq:choice-beta} above and rearranging to solve for $p$ gives
\begin{equation}
\label{eq:choice-p}
p = \left( \frac{\beta(\beta+1)c_\sigma \sigma}{\pi} \right)^2.
\end{equation}
As in the rounding algorithm, we require in addition that $p>c_1k$, so the sample lengths for the multiscale algorithm are chosen to satisfy
\begin{equation}\label{eq:multi-p}
 p > \max\left\{ c_1k,  \left( \frac{\beta(\beta+1)c_\sigma \sigma}{\pi} \right)^2 \right\}. 
\end{equation}

\subsubsection{Number of iterations}
\label{subsubsec:num-it}
Recall from Corollary~\ref{coro-4.3} that, for constant $\beta_j = \beta$, $m= \left\lfloor \log_\beta \frac{2\delta}{p\eps_0} \right\rfloor+1$ shifts suffices to ensure that the estimated frequency satisfies $|\wt\omega-\omega|<\frac{p}{2}$. As in Section~\ref{sec:minor} we take $\eps_0 = \frac{1}{2N}$ to avoid the branch of $\Arg(z)$. Assume that the first term in~\eqref{eq:multi-p} is the larger of the two, so that $p=\bigo(k)$. Then after $\bigo(\log(N/k))$ iterations, by rounding the approximate frequency $\wt\omega$ to the closest integer of the form $np+h$, where $h=\omega\wmod{p}$ is known from the location of the peak in $\wh{\bS}_p^\rmn$, we will recover the true frequency $\omega$. With the results of \cite{lawlor2013adaptive} this immediately implies that the average-case runtime of the multiscale algorithm is $\bigo(k \log(k) \log(N/k))$.

\subsubsection{Robust aliasing test}
\label{subsubsec:robust-aliasing}
As noted in Subsection~\ref{subsec:notation}, our frequency estimation procedure works only for non-collsion $\omega$. In \cite{lawlor2013adaptive} two tests were given to determine whether a collision had occurred at a candidate frequency. In the implementation of that algorithm in the noiseless setting, requiring the ratio~\eqref{eq:non-collision} to be within some threshold of unity sufficed to detect collisions. In the setting of the current paper, where the samples are corrupted with noise, we resort to the second of the tests given in in \cite{lawlor2013adaptive}, which examines the ratios~\eqref{eq:non-collision} for several values of $\eps$. For $0\le j \le m$ we compute the ratio~\eqref{eq:non-collision} and compare it with a threshold $\tau$. We count the number of times the ratio exceeds $\tau$  and reject those frequencies which fail more than an $\eta$ fraction of the tests. Since we expect fluctuations in this ratio due to noise of order $\sigma/\sqrt{p}$ we set $\tau$ to be a small constant multiple of this quantity. 

We give pseudocode for the iterative frequency estimation procedure below; the full
algorithm is given by replacing the calculation of frequencies in \cite[Algorithm 1]{lawlor2013adaptive}
with this procedure.

%\begin{algorithm}
%    \caption{\textsc{MultiscaleFreqEst}}
%    \label{alg:est}
    \begin{algorithmic}[5]
        \REQUIRE $S(t), N, k, \beta, \sigma, c_\sigma, \eta$
        \ENSURE $\{ \wt{\omega}_\ell \}_{\ell=1}^{\wt{k}}$
        \STATE $p \gets \max\left\{ c_1k,  \left( \frac{\beta(\beta+1)c_\sigma \sigma}{\pi} \right)^2 \right\}$
        \STATE $\tau \gets \frac{c_\sigma \sigma}{\sqrt{p}}, m \gets 1 + \left\lfloor \log_\beta \frac{N}{p}\right\rfloor$
        \STATE $\textup{vote}_\ell \gets 0, \; \ell=1, \ldots, k$
        \STATE $\wh{\bS}_{p} \gets$ FFT of $\frac{1}{p}$-samples of $S(t)$

        \FOR{$j=0$ to $m$}

            \STATE $\eps_j \gets \frac{\beta^j}{2N}$
            \STATE $\wh{\bS}_{p,\eps_j} \gets$ FFT of $\eps_j$-shifted $\frac{1}{p}$-samples of $S(t)$

            \FOR {$\ell=1$ to $k$}
                \STATE $h \gets$ index of $\ell^\textup{th}$ largest peak in $\wh{\bS}_p$
                \STATE $r \gets \left| \frac{\left|\wh{\bS}_{p,\eps_j}[h] \right|}{\left|\wh{\bS}_{p}[h] \right|} -1 \right|$
                \IF {$r>\tau$}
                    \STATE $\textup{vote}_\ell \gets \textup{vote}_\ell +1$
                \ENDIF
                \STATE $b_j \gets \frac{1}{2\pi}\Arg\left( \frac{\wh{\bS}^\mathrm{n}_{p,\eps_j}[h]}{\wh{\bS}^\mathrm{n}_p[h]} \right)$
                \IF {$j=0$}
                    \STATE $\wt{\omega}_\ell \gets b_j/\eps_j$
                \ELSE
                    \STATE $\wt{\omega}_\ell \gets \wt{\omega}_\ell + (b_j-\eps_j\wt{\omega}_\ell)\wmodc/\eps_j$
                \ENDIF
                \IF {$j=m$}
                    \STATE $\wt{\omega}_\ell \gets p\cdot\textup{round}\left( \frac{\wt{\omega}_\ell-h}{p} \right) +h$
                \ENDIF
            \ENDFOR
        \ENDFOR
        \STATE return $\wt{\omega}_\ell$ with $\textup{vote}_\ell \le \eta(m+1)$
    \end{algorithmic}
%\end{algorithm}

\section{Empirical evaluation}
\label{sec:empirical}

In this section we describe the results of an empirical evaluation of the algorithms of section~\ref{sec:minor} and \ref{sec:multi}. We focus on two aspects of the algorithms' performance: accuracy as measured in the EMD(1) and EMD($\omega$) metrics (c.f.~Subsection~\ref{subsec:emd}), and runtime as a function of both the sparsity $k$ and the noise level $\sigma$. In all of the experiments reported below, we report averages over 100 random test signals generated according to the prescription in Subsection~\ref{subsec:random_signal_model}. The bandwidth for these tests was fixed at $N=2^{22}$.

All experiments were conducted in C++ on a Linux machine with four Intel Xeon X5355 dual-core processors at 2.66 GHz and 64 Gb of RAM. The GNU compiler was used with optimization flag -O3. For the multiscale algorithm, it was determined after extensive testing that the choice of parameters $c_1=2, \; c_\sigma=6, \; \eta=\frac{1}{4}, \beta=2.5$ gave a satisfactory balance between runtime and accuracy. All FFTs are performed using FFTW3 \cite{frigo2005design}. For comparison, we also present the results of the same trials for two alternative sparse Fourier algorithms: sFFT 1.0 \cite{hassanieh2012simple} and AAFFT \cite{iwen2007empirical}. 

\subsection{Accuracy}
\label{subsec:accuracy}
\begin{figure}
    \centering
    \subfloat[]{\includegraphics[width=0.9\columnwidth]{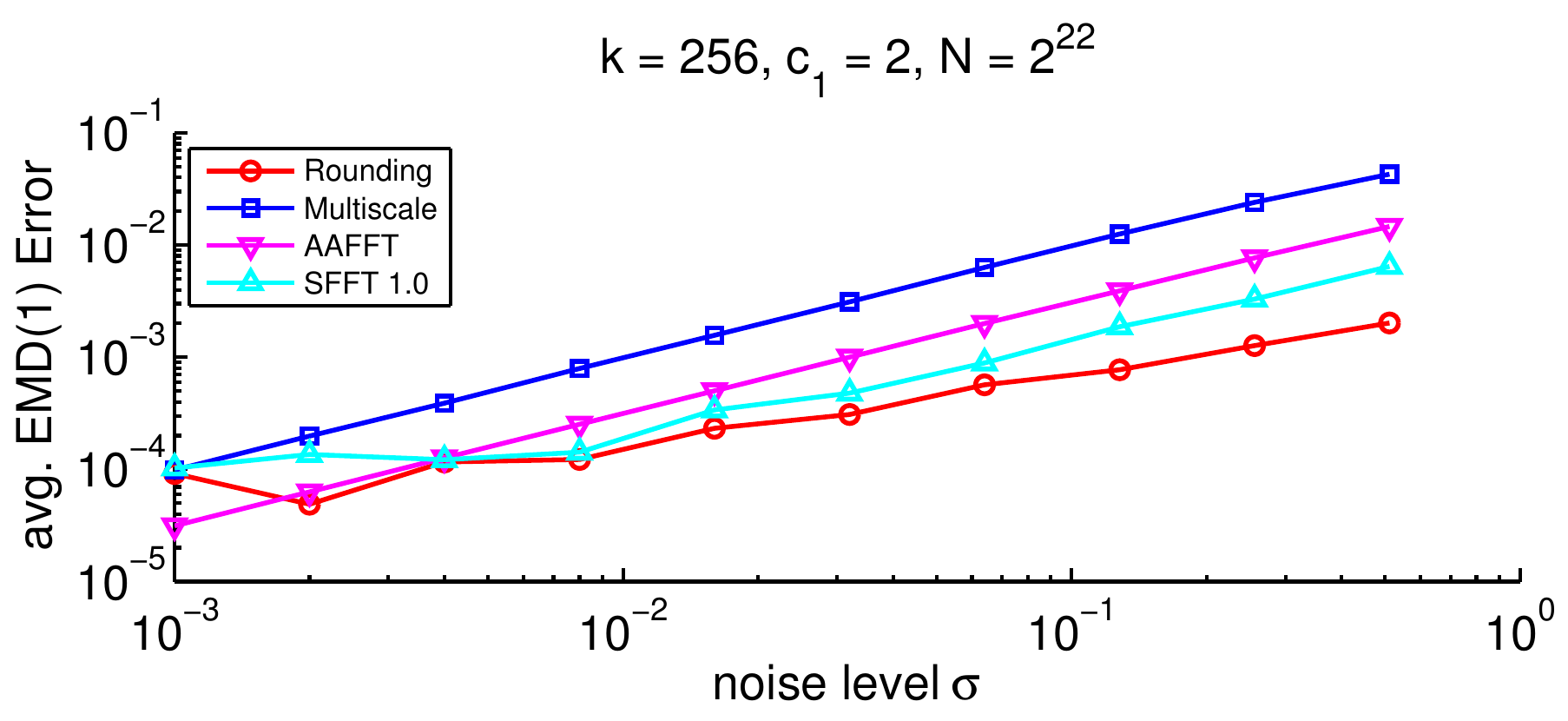}}\\
    \subfloat[]{\includegraphics[width=0.9\columnwidth]{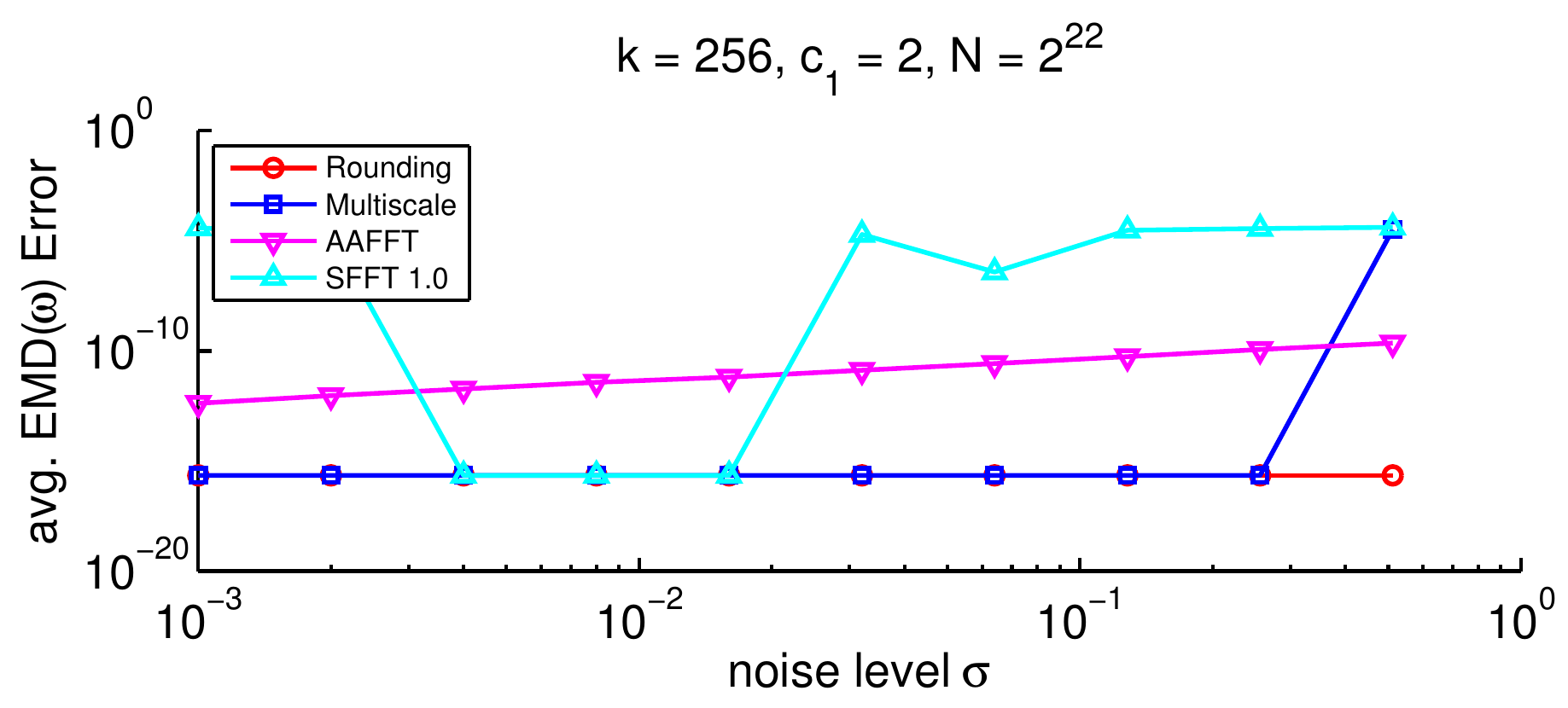}}
    \caption{(a) Average EMD(1) error of the algorithms as a function of noise level $\sigma$. (b) Average EMD($\omega$) error as a function of noise level $\sigma$. Due to the log scale on the $y$ axis, all EMD($\omega$) values have been shifted up by $10^{-16}$ for clarity.}
    \label{fig:error}
\end{figure}

In Figure~\ref{fig:error} (a) we plot the average EMD(1) error of the algorithms as a function of the noise level $\sigma$. For the rounding algorithm, the EMD(1) error increases as $\sigma^{2/3}$, while for the other three it increases linearly.  In all cases the EMD(1) error is dominated by the coefficient error. The coefficient estimates in all four algorithms are given by an empirical average of the samples, and so the accuracy is determined by the number of samples taken. This explains both the scaling of the error of our rounding algorithm (recall from Section~\ref{sec:minor} that $p>(\sigma/\eps)^{2/3}$), as well as the larger EMD(1) error of our multiscale algorithm, which performs well even with $c_1$ as small as two. The multiscale error correction allows us to take much coarser sampling rates to achieve a tolerable error. As we show in the next subsection, these coarser sampling rates lead to much improved runtime.

In order to assess the accuracy of the frequency lists returned by each of the four algorithms, in Figure~\ref{fig:error} (b) we plot the average EMD($\omega$) error as a function of the noise level. The EMD($\omega$) error was zero for all trials of the rounding algorithm, as expected due to the choice of $p$. Moreover, for all but the highest noise level, the EMD($\omega$) error of the multiscale algorithm was zero in all trials. For most values of $\sigma$, the EMD($\omega$) error of sFFT 1.0 was non-zero, indicating that even at low to moderate noise levels, erroneous frequencies are returned. The EMD($\omega$) error of AAFFT was always less than $1/N$, indicating that true frequencies were recovered in all cases; the non-zero values are numerical artifacts.

\subsection{Runtime}
\label{subsec:runtime}
\begin{figure}
    \centering
    \subfloat[]{\includegraphics[width=0.9\columnwidth]{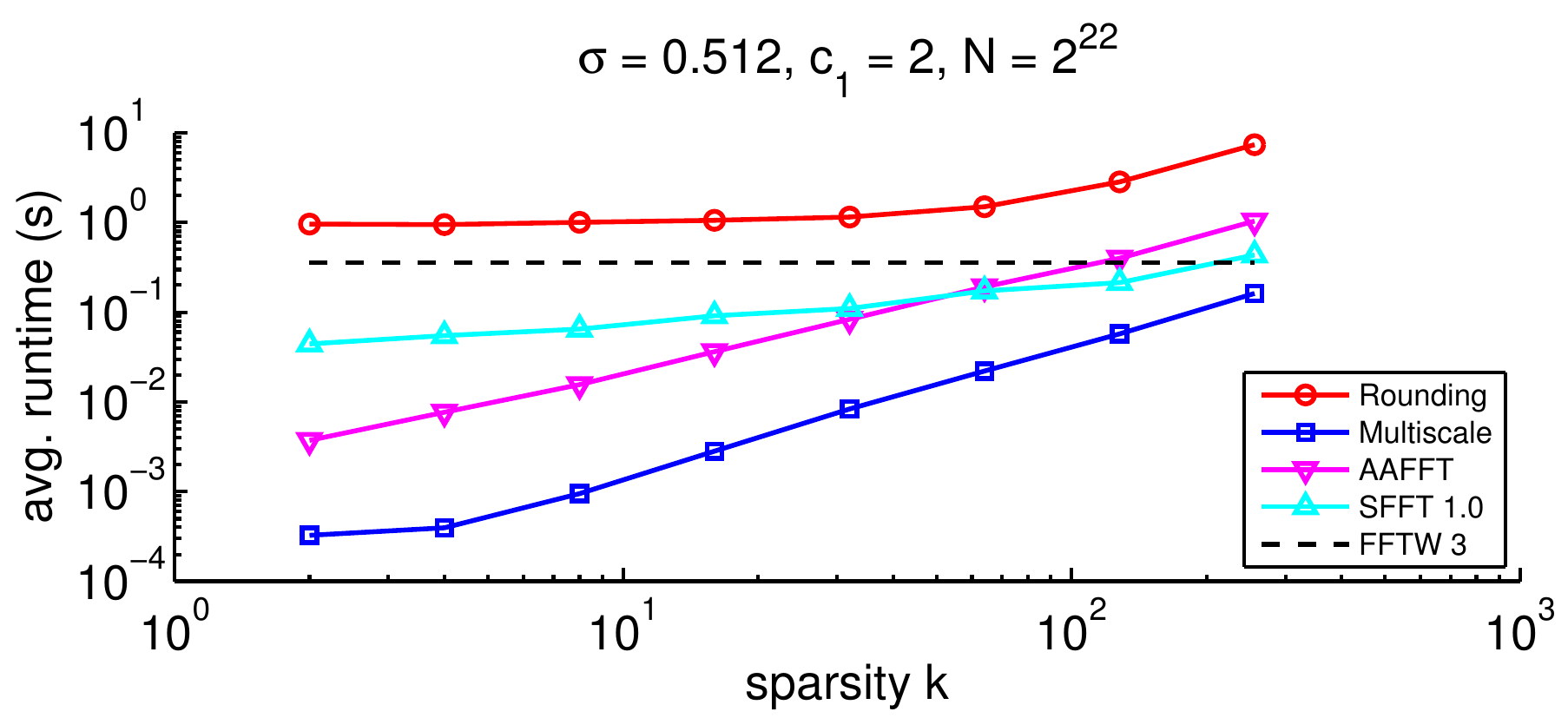}}\\
    \subfloat[]{\includegraphics[width=0.9\columnwidth]{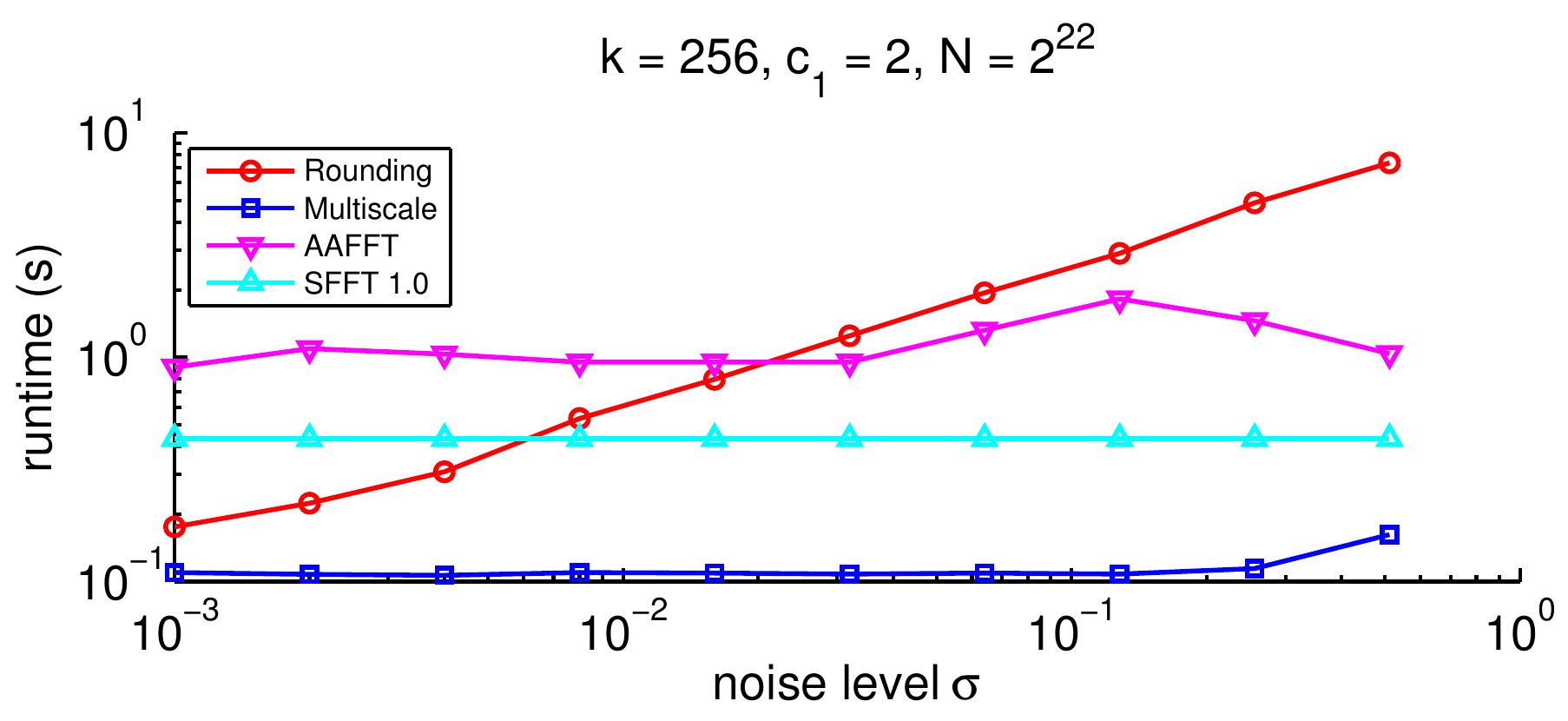}}
    \caption{(a) Average runtime vs. sparsity $k$ for the algorithms tested. (b) Average runtime vs. noise level $\sigma$.}
    \label{fig:runtime}
\end{figure}

In Figure~\ref{fig:runtime} (a) we plot the average runtime of the algorithms as a function of the sparsity $k$ for a fixed value of the noise level $\sigma=0.512$ and the parameter $c_1=2$. As a reference for runtime comparisons, we also plot the time taken by FFTW3 on the same machine. For the rounding algorithm, we see that there is no dependence on $k$ until $k = 64$; this is a consequence of the requirement \eqref{eq:new_p} on the choice of sampling rate.  Thus at this noise level our modified algorithm is slightly slower than a highly optimized FFT implementation.  The average runtime of our multiscale algorithm scales slightly superlinearly with $k$, which is expected given the runtime bound $\bigo(k \log (k) \log(N/k))$ of Subsection \ref{subsubsec:num-it}. Moreover, we note that for all levels of sparsity tested, the multiscale algorithm outperforms AAFFT, sFFT 1.0, and FFTW3.

In Figure \ref{fig:runtime} (b) we plot the average runtime of the algorithms as a function of the noise level $\sigma$ for a fixed value of the sparsity $k=256$. For the rounding algorithm we can see the approximate dependence of the runtime on $\sigma^{2/3}$, as dictated by the choice of $p$ in \eqref{eq:new_p}. For the multiscale algorithm, there is no dependence on $\sigma$ until the very noisy case $\sigma=.512$. 

\subsection{Spurious frequencies}
\begin{figure}
    \centering
    \includegraphics[width=0.9\columnwidth]{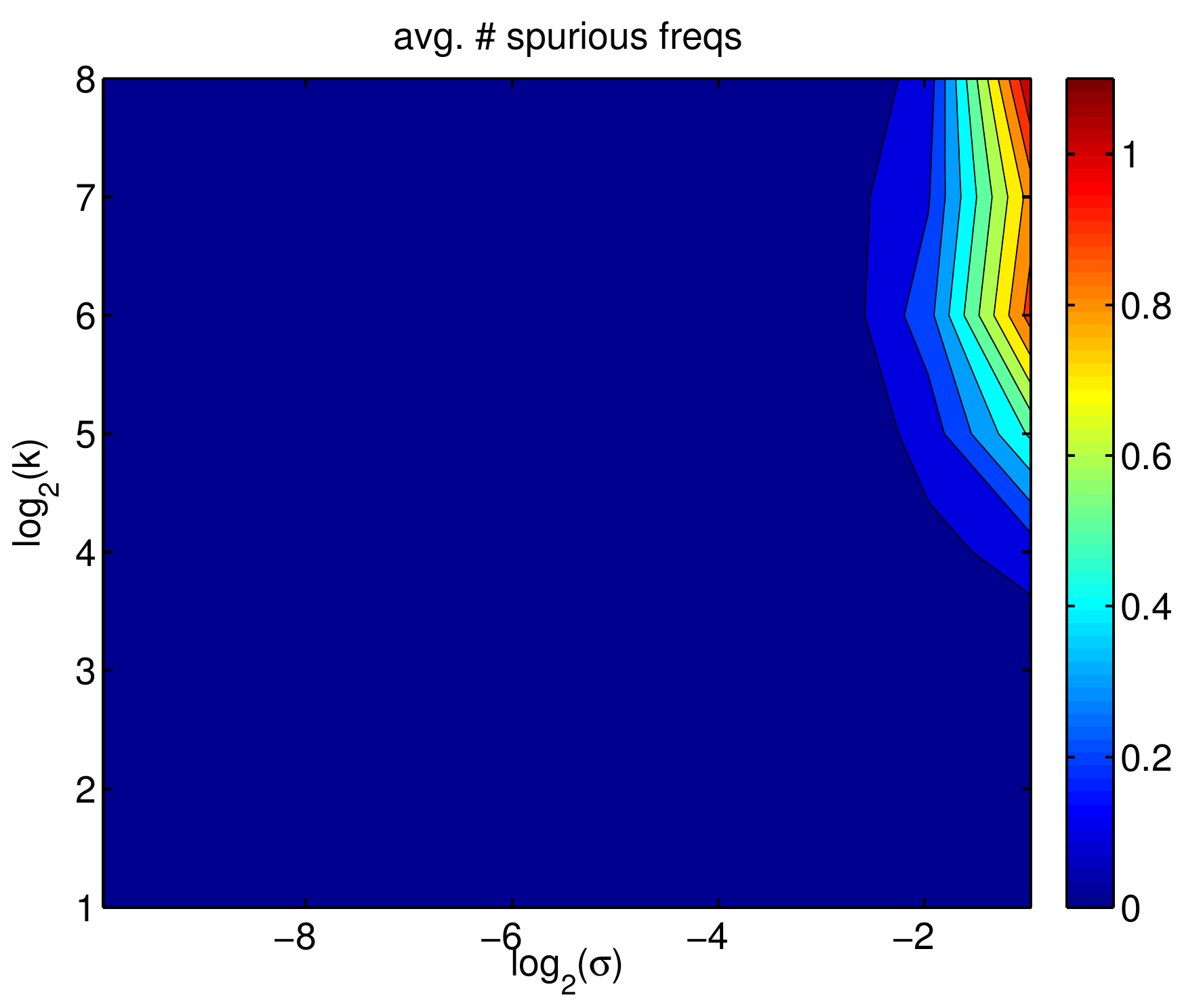}
    \caption{Average number of spurious frequencies inserted and deleted by the multiscale algorithm as a function of $k$ and $\sigma$.}
    \label{fig:spurious}
\end{figure}

As noted in Section~\ref{sec:prelim}, due to noise it is possible that one or more spurious frequencies are introduced into our signal representation. In subsequent iterations, it is extremely likely that any such spurious frequency will be identified and subtracted from the updated representation. Since this happens with non-zero probability, it is of interest to examine how often such an insertion and deletion occurs. In Figure~\ref{fig:spurious}, we plot the average number of spurious frequencies inserted and deleted by the multiscale algorithm as a function of $k$ and $\sigma$. It is clear that the inclusion of spurious frequencies only occurs in the high-noise, high-sparsity regime. Moreover, on average only one such wrong frequency appears in our representation even in this challenging environment. This indicates that our robust aliasing test of Subsection~\ref{subsubsec:robust-aliasing} does a very good job at detecting collisions in all but the most extreme circumstances.

%%%%%%%%%%%%%%%%%%%%%
\section{Conclusion}
\label{sec:conc}
In this paper we gave two extensions of the sparse Fourier algorithm of \cite{lawlor2013adaptive} to handle noisy signals. The first of these was a minor modification of the original algorithm that involved rounding frequency estimates to the nearest integer with the correct residue modulo the sampling rate. We showed that in order for this modification to correctly identify the true frequencies in Gaussian noise of standard deviation $\sigma$ the sampling rate needed to satisfy $p \ge \sigma^{2/3}$. While this resulted in accurate approximations of the Fourier transform in the EMD(1) and EMD($\omega$) metrics, the sampling rate requirement forced the algorithms to be slow in practice.

The second extension overcame this pitfall by introducing a novel multiscale approach to frequency estimation in the sparse Fourier transform context. By using samples of the input at multiple time shifts spaced geometrically, our algorithm exhibits a form of error correction in its frequency estimation. This allows the use of much coarser sampling rates than the first modification, which in turn leads to greatly reduced runtimes in our empirical evaluation. The error correction of our multiscale algorithm is to the best of our knowledge novel in the sparse Fourier transform context, and we believe it is a promising approach for further investigation.

%%%%%%%%%%%%%%%%%%%%%%%%%%
\section*{Acknowledgments}
During the preparation of this manuscript we became aware of related work by Laurent Demanet and his collaborators. We kindly acknowledge their generosity in sharing their work with us. We also thank the anonymous reviewers whose suggestions improved the exposition of this manuscript.

%%%%%%%%%%%%%%%%%%%%%%%%%%%
\bibliographystyle{amsalpha}
\bibliography{papertwo}

%%%%%%%%%%%%%%
\end{document}